\documentclass[12pt,centertags,oneside]{article}
\usepackage{amsmath,amstext,amsthm,amscd,typearea}
\usepackage{amssymb}
\usepackage{a4wide}
\usepackage[mathscr]{eucal}
\usepackage{mathrsfs}
\usepackage{typearea}
\usepackage{charter}
\usepackage{pdfsync}
\usepackage{url}

\usepackage{xcolor}

\usepackage[a4paper,width=16.2cm,top=3cm,bottom=3cm]{geometry}

\numberwithin{equation}{section}

\newtheorem{theorem}{Theorem}[section]

\newtheorem{proposition}[theorem]{Proposition}
\newtheorem{corollary}[theorem]{Corollary}
\newtheorem{lemma}[theorem]{Lemma}

\newcommand{\cali}[1]{\mathscr{#1}}

\newcommand{\vol}{\mathop{\mathrm{vol}}}

\newcommand{\ddc}{dd^c}
\newcommand{\dc}{d^c}

\newcommand{\PSH}{{\rm PSH}}

\newcommand{\capK}{\text{cap}}

\newcommand{\B}{\mathbb{B}}

\newcommand{\N}{\mathbb{N}}

\newcommand{\R}{\mathbb{R}}

\title{\bf  Complex Monge-Amp\`ere equations with solutions in finite energy classes}
\providecommand{\keywords}[1]{\textbf{\textit{Keywords:}} #1}
\providecommand{\subject}[1]{\textbf{\textit{Mathematics Subject Classification 2010:}} #1}

\author{Do Duc Thai and Duc-Viet Vu}

\newcommand{\Addresses}{{
		\bigskip
		\footnotesize
		\textsc{Duc-Viet Vu, University of Cologne, Division of Mathematics, Department of Mathematics and Computer Science, Weyertal 86-90, 50931, K\"oln.}
		\noindent
		\par\nopagebreak
		\noindent
		\textit{E-mail address}: \texttt{vuviet@math.uni-koeln.de}
\newline
		
		\textsc{Do Duc Thai, Department of Mathematics, Hanoi National University of Education, 136 XuanThuy str., Hanoi, Vietnam.}
		\noindent
		\par\nopagebreak
		\noindent
		\textit{E-mail address}: \texttt{doducthai@hnue.edu.vn}	
}}

\date{\today}
\begin{document}
\maketitle
\begin{abstract} We characterize the class of probability measures on  a compact K\"ahler manifold such that the associated Monge-Amp\`ere equation has a solution of finite pluricomplex energy. Our results are also valid in the big cohomology class setting.
\end{abstract}
\noindent
, {pluricomplex energy}, {full mass intersection}.
\\

\noindent
\subject{32U15},  {32Q15}.



\section{Introduction}

 Let $X$ be a compact K\"ahler manifold of dimension $n$ and $\omega$ a K\"ahler form on $X$ normalized so that $\int_X \omega^n =1$.  
We are interested in studying the regularity of the solutions of the complex Monge-Amp\`ere equation
\begin{align}\label{eq-MA}
(\ddc \varphi+ \omega)^n =\mu,
\end{align}
where $\mu$ is a probability measure on $X$, and $\varphi$ is an $\omega$-psh function. The case where $\mu$ is smooth was settled in \cite{Yau1978}. The pluripotential theory tools were introduced to the study of (\ref{eq-MA}) in \cite{Kolodziej_Acta}, and  an almost sharp class of measures whose equation (\ref{eq-MA}) admits a continuous solution was established there; see also \cite{NgocCuong-Kolodziej,KC_remark_continuous}. 

The class of measures with H\"older continuous solutions for  (\ref{eq-MA}) was characterized by contributions of many works. Roughly speaking, the equation (\ref{eq-MA}) has a H\"older continuous solution if and only if certain functional associated to $\mu$ is H\"older continuous. Direct analogue of this fact in the continuity regularity is not true. We refer to \cite{DemaillyHiep_etal,DKC_Holder-Sobolev,DinhVietanhMongeampere,KC-holder2018,Kolodziej08holder,Lu-To-Phung,NgocCuong-Holder-calc.var,NgocCuong-Holder2020,Tosatti-Weinkove,Vu_MA,Vu_MA_holder} and references therein for details. In this paper, we are interested in a weaker notion of  regularity for solutions of  (\ref{eq-MA}) which is the so-called finite pluricomplex energy. The last term was coined in \cite{Cegrell}.  Precisely, our goal is to characterize  the class of measures $\mu$ such that  (\ref{eq-MA}) has solutions of finite pluricomplex energy.  To go into details, we need to recall some notions. 

For closed positive $(1,1)$-currents $T_1, \ldots, T_m$, let $\langle T_1 \wedge \cdots \wedge T_m \rangle $ be  the \emph{non-pluripolar product} of $T_1, \ldots, T_m$. 
We refer to \cite{BT_fine_87,BEGZ,Lu-Darvas-DiNezza-mono,GZ-weighted,Viet-convexity-weightedclass,Viet-generalized-nonpluri,WittNystrom-mono,Xia} and references therein for basic properties of  non-pluripolar products. 

Let $\mathcal{E}$ be the set of negative $\omega$-psh functions $\varphi$ of full Monge-Amp\`ere mass, \emph{i.e.}, $\big\langle (\ddc \varphi+ \omega)^n\big \rangle$ is of mass equal to $1$. Let $\mathcal{W}^-$ be the set of convex increasing functions $\chi: \R \to \R$ such that $\chi(0)=0$ and  $\chi(-\infty)= -\infty$. Let $M\ge 1$ be a constant and $\mathcal{W}^+_M$ be the set of  concave increasing  functions $\chi: \R \to\R$ such that $\chi(0)=0$,  $\chi(-\infty)= -\infty$ and $|t\chi'(t)| \le M |\chi(t)|$ for $t \in \R^{-}$. 

For $\chi \in \mathcal{W}^- \cup \mathcal{W}^+_M$, we denote by $\mathcal{E}_\chi$ the space of  $\omega$-psh functions $\varphi \in \mathcal{E}$ such that 
$$E_\chi(\varphi):= - \int_X \chi(\varphi) \big \langle (\ddc \varphi+ \omega)^n \big \rangle< \infty.$$
The last quantity is called \emph{the $\chi$-energy} of $\varphi$. The space  $\mathcal{E}_\chi$ was introduced in \cite{GZ-weighted} as a counterpart in the compact setting of  the class of psh functions of finite pluricomplex energy given in \cite{Cegrell}. For a positive real number $p$ and $\chi(t)= -(-t)^p$, we write $\mathcal{E}^p,E_p$ for $\mathcal{E}_\chi, E_\chi$ respectively.  Recall that 
\begin{align} \label{eq-EEchi}
\mathcal{E}= \cup_{\chi \in \mathcal{W}^-} \mathcal{E}_\chi, \quad \mathcal{E}_{\chi_1} \subset \mathcal{E}^{1} \subset \mathcal{E}_{\chi_2},
\end{align}
for every $\chi_1 \in \mathcal{W}^+_M$ and $\chi_2 \in \mathcal{W}^-$; see \cite{GZ-weighted}.  

By \cite{Cegrell,GZ-weighted},  for every probability measure $\mu$ having no mass on pluripolar sets, the equation (\ref{eq-MA}) has a solution $\varphi \in \mathcal{E}$ which is in fact unique by \cite{Dinew-uniqueness}. We refer to \cite{AhagCegrellHiepMA} for more information in the local setting.  Hence, it is natural to characterize the class of measures $\big \langle (\ddc \varphi+ \omega)^n \big \rangle$ for $\varphi \in \mathcal{E}_\chi$; see the comments after Theorem \ref{th-main} for details. When $\chi(t)= -(-t)^p$ for some constant $p>0$ (or $\chi$ is a quasi-homogeneous weight), the question was solved in \cite{Cegrell,GZ-weighted} thanks to the special form of $\chi$.  Here is our main result giving a complete answer for the general case.

\begin{theorem} \label{th-main}  Let $\mu$ be a probability measure on $X$. Let $\chi \in \mathcal{W}^- \cup \mathcal{W}_M^+$. Then, $$\mu= \big \langle (\ddc \varphi+ \omega)^n\big \rangle$$ for some $\varphi \in \mathcal{E}_\chi$ if and only if there is a constant $A>0$ such that 
\begin{align}\label{ine-dktuongduongEchi}
- \int_X \chi(\psi) d \mu \le  A\big[\big(E_\chi(\psi)\big)^\lambda+1\big]
\end{align}
for every $\psi \in \mathcal{E}_\chi$ with $\sup_X \psi=-1$, where $\lambda:= 1/2$ if $\chi \in \mathcal{W}^-$, and $\lambda:= M/(M+1)$ if $\chi \in \mathcal{W}^+_M$. In particular, the set of measures $\big \langle (\ddc \varphi+ \omega)^n \big \rangle$ with $\varphi \in \mathcal{E}_\chi$ is convex.  
\end{theorem}

We note that Theorem \ref{th-main} is closely related to a problem posed  in \cite{GZ-weighted} (see the comment right after \cite[Theorem 4.1]{GZ-weighted}): for every probability measure $\mu$, does the the condition that $\chi(\psi) \in L^1(\mu)$ for every  $\psi \in \mathcal{E}_\chi$ imply that $\mu= \big \langle (\ddc \varphi+ \omega)^n\big \rangle$ for some $\varphi \in \mathcal{E}_\chi$? Although  Theorem \ref{th-main} doesn't answer that question, the important point is that it gives a characterization of the class of measures $\big \langle (\ddc \varphi+ \omega)^n \big \rangle$ for $\varphi \in \mathcal{E}_\chi$ by an integrability property. Such a condition is of great interest in determining the range of complex Monge-Amp\`ere operators. We refer to \cite[Question 13-15]{DGZ} for a discussion regarding this problem. 
 
Previously,  it was shown in \cite[Theorem 3.1]{Benelkourchi-Guedj-Zeriahi} that if $\mu$ satisfies certain conditions in terms of capacity, then (\ref{eq-MA}) has a solution in finite energy classes.  In the local setting of hyperconvex domains,  a result similar to Theorem \ref{th-main}  for $\chi \in \mathcal{W}^-$ was obtained in \cite[Theorem C]{Benelkourchi}.  
 

By \cite[Proposition 3.2]{BEGZ}, for every non-pluripolar probability measure $\mu$ on $X$, there exists $\chi \in \mathcal{W}^-$ such that $- \int_X \chi(\psi) d \mu< \infty$ for every $\omega$-psh function $\psi$ with $\sup_X \psi =-1$, in particular,  (\ref{ine-dktuongduongEchi}) is satisfied. Using this, (\ref{eq-EEchi}) and Theorem \ref{th-main}, we recover a well-known fact proved in \cite{GZ-weighted} that every non-pluripolar probability measure can be written as $\big \langle (\ddc \varphi+ \omega)^n\big \rangle$ for some $\varphi \in \mathcal{E}$. This is an instance showing why we are interested in the class $\mathcal{E}_\chi$ for $\chi \in\mathcal{W}^-$. Broadly speaking, we can study non-pluripolar measures by considering them as Monge-Amp\`ere measures with potentials  in some class $\mathcal{E}_\chi$ for $\chi \in \mathcal{W}^-$ which is easier to handle analytically. This point of view should be useful in complex dynamics because the equilibrium measure associated to a meromorphic self-map is usually non-pluripolar; see, for example, \cite{Diller-Guedj-Dujardin,DS_acta,GZ-weighted,Vu_nonkahler_topo_degree}.   In this aspect, considering only $\chi(t)= -(-t)^p$ for $0<p\le 1$ is not enough  because in general the union of  $\mathcal{E}_\chi$ with $\chi(t)= -(-t)^p$  ($0<p\le 1$) is a proper subset of $\mathcal{E}$ by \cite[Example 2.13]{GZ-weighted}. 

Another reason why it is of interest to study general weights $\chi$  rather than only those $\chi(-t)=-(-t)^p$ for $p>0$ is that in some problems even when one only considers the weight $-(-t)^p$, it is still necessary to approximate this weight by those in $\mathcal{W}^- \cup \mathcal{W}_M^+$. To illustrate this, one can see the proofs of \cite[Theorem 2.17]{BEGZ}, \cite[Theorem 3.5]{Darvas-finite-energy} and the paragraph after \cite[Theorem 1.1]{Darvas-Lu-Rubinstein} in \cite{Darvas-Lu-Rubinstein}.

As a direct consequence of Theorem \ref{th-main}, we obtain immediately the following monotonicity property. 

\begin{corollary} \label{cor-subsolutionkahler} Let $\chi \in \mathcal{W}^- \cup \mathcal{W}_M^+$.  Let $\mu$ be a probability measure such that 
$$\mu \le  C  \big \langle(\ddc \varphi+ \omega)^n \big \rangle$$
 for some $\varphi \in \mathcal{E}_\chi$ and some constant $C$. Then there exists $\psi \in \mathcal{E}_\chi$ satisfying 
$$\mu= \big \langle (\ddc \psi+ \omega)^n  \big \rangle.$$ 
\end{corollary}

We note that the above result fits in the study of a long-standing conjecture due to Ko{\l}odziej that if $\mu$ is bounded by a Monge-Amp\`ere measure with bounded or continuous potentials then $\mu$ is also a Monge-Amp\`ere measure with bounded or continuous potentials (\cite[Question 15]{DGZ},  \cite[Theorem 4.7]{Kolodziej05} and \cite{KC_remark_continuous}). It reflects a more general belief that if a closed positive current $T$ is dominated by another current $S$, then $T$ should inherit properties from $S$. Another example of this intuition is \cite[Theorem 1.1]{DNV} proving that if $S$ has H\"older, continuous, or bounded super-potentials, then so does $T$.  

We have some comments on the proof of Theorem \ref{th-main}. We follow the standard strategy in \cite{Cegrell,GZ-weighted,Kolodziej_Acta}. The idea goes as follows. Firstly, we regularize $\mu$ to obtain a sequence of measures $(\mu_j)_j$ with nicer regularity. We solve the Monge-Amp\`ere equation for $\mu_j$ to get a solution $\varphi_j$. The $L^1$ limit of the sequence $(\varphi_j)_j$ is supposed to be the solution we are looking for.   A key point in this approach is to control the regularity of the sequence of $(\varphi_j)_j$. In our present setting, this means that we need to prove that the $\chi$-energy of $\varphi_j$ is bounded uniformly in $j$. This was achieved for $\chi(t)= -(-t)^p$ in previous works. Our contribution is to prove that this can be done for general $\chi$. The key ingredient is Theorem \ref{th-integra} below giving an upper bound for ``mixed" energy. 

Our method also works in the big cohomology class setting and gives the following generalization of Theorem \ref{th-main}. 

\begin{theorem} \label{th-mainbigclass} Let $\theta$ be a smooth real $(1,1)$-form in a big cohomology class in $X$.  Let $\phi$ be a $\theta$-psh function  of model type singularity on $X$ with $\int_X \big \langle  (\ddc \phi+\theta)^n  \big \rangle >0$. Let $\chi \in \mathcal{W}^- \cup \mathcal{W}_M^+$. Then, 
$$\mu= \big \langle (\ddc \varphi+ \theta)^n \big \rangle$$
 for some $\varphi \in \mathcal{E}_\chi(\theta,\phi)$ if and only if there is a constant $A>0$ such that 
\begin{align}\label{ine-dktuongduongEchibig}
- \int_X \chi(\psi-\phi) d \mu \le   A\big[\big(E_{\chi,\theta,\phi}(\psi)\big)^\lambda+1\big] 
\end{align}
for every $\psi \in \mathcal{E}_\chi(\theta,\phi)$ with $\sup_X (\psi-\phi)=-1$, where $\lambda:= 1/2$ if $\chi \in \mathcal{W}^-$, and $\lambda:= M/(M+1)$ if $\chi \in \mathcal{W}^+_M$. 
\end{theorem}

For the definition of model type singularity and other notations in the statement of Theorem \ref{th-mainbigclass}, we refer to Section \ref{sec-bigcoho}. The condition $\int_X \big \langle (\ddc \phi+\theta)^n \big \rangle >0$ is natural as explained in \cite{Lu-Darvas-DiNezza-logconcave,Lu-Darvas-DiNezza-mono} (note that if $\int_X \big \langle (\ddc \phi+\theta)^n\big \rangle=0$, then by monotonicity, $\big \langle (\ddc \varphi+ \theta)^n \big \rangle=0$ for every $\varphi \in \mathcal{E}_\chi(\theta, \phi)$). The hypothesis imposed on $\phi$ is always satisfied if $\phi$ is of minimal singularities in the cohomology class of $\theta$. 

We note here that by \cite{Lu-Darvas-DiNezza-logconcave,Lu-Darvas-DiNezza-mono}, every measure having no mass pluripolar sets can be written as $ \big \langle (\ddc \varphi+ \theta)^n \big \rangle$ for some $\varphi \in \mathcal{E}(\theta,\phi)$. However, it seems that the arguments there are not enough to obtain the above characterization given in Theorem \ref{th-mainbigclass} even when $\chi(t)= -(-t)^p$. 

In order to be able to transfer the arguments from the K\"ahler case to the big case,  it is crucial to have an adequate integration by parts formula. Such a formula was proved  in \cite[Theorem 1.14]{BEGZ} under the hypothesis of having small unbounded locus. There are also some partial generalizations given in \cite{Lu-cap-compare,Xia}.  However, they are not good enough to treat  the setting of prescribed singularity type as in Theorem \ref{th-mainbigclass}. To overcome this problem, we use an integration by parts formula obtained \cite{Viet-convexity-weightedclass}  (see Theorem \ref{th-integrabypart} below) in an essential way. This formula fully extends \cite[Theorem 1.14]{BEGZ} to the context of prescribed singularity type.  As another illustration of the use of Theorem \ref{th-integrabypart}, we will explain how to use the variational method in \cite{BBGZ-variational} to  give another proof  of \cite[Theorem A]{Lu-Darvas-DiNezza-logconcave} which was proved there using supersolution approach; see Theorem \ref{th-variation} below.

Finally, one can see that Theorem \ref{th-mainbigclass} implies a generalization of Corollary \ref{cor-subsolutionkahler} in the big cohomology class setting. Theorem \ref{th-main} is proved in the next section.  We prove Theorem \ref{th-mainbigclass} in the last section. \\

\noindent
\textbf{Acknowledgments.} We thank the referees for suggestions improving considerably the presentation of the paper. D.-V. Vu is supported by  a postdoctoral fellowship of the Alexander von Humboldt Foundation.

\section{Proof of Theorem \ref{th-main}}


If no confusion can arise, we sometimes use $\gtrsim, \lesssim$ to denote $\ge, \le$ modulo multiplicative constants independent of parameters in question. We begin with some elementary lemmas, which appeared more or less in \cite[Section 2.1]{Darvas-finite-energy}. 

\begin{lemma} \label{le-regularizedchi} Let $\chi \in \mathcal{W}^- \cup \mathcal{W}_M^+$.  Let $g$ be a smooth radial cut-off function supported in $[-1,1]$ on $\R$, \emph{i.e,} $g(t)= g(-t)$ for $t \in \R$, $0 \le g \le 1$ and $\int_\R g(t) dt =1$. Put $g_\epsilon(t):=  \epsilon^{-1}g(\epsilon t)$ for every constant $\epsilon >0$ and $\chi_\epsilon:= \chi * g_\epsilon$ (the convolution of $\chi$ with $g_\epsilon$). Then $\chi_\epsilon$ converges uniformly to $\chi$ as $\epsilon \to 0$ on compact subsets in $\R$, and  the following assertions are true:

$(i)$ if $\chi \in \mathcal{W}^-$, then   $\chi_\epsilon- \chi_\epsilon(0) \in \mathcal{W}^-$ and $\chi_\epsilon \ge \chi$,
 
 $(ii)$ if $\chi \in \mathcal{W}_M^+$, then  $\chi_\epsilon$ is concave increasing,  $\chi_\epsilon \le \chi$ and  for every constant $\delta>0$, we have
 $$|t \chi'_\epsilon(t)|  \le -\big(M+ \epsilon/ (\delta- \epsilon)\big) \chi_\epsilon(t)$$
 for $\epsilon< \delta$ and  $t \le -\delta$. 
\end{lemma}

\proof  By the continuity of $\chi$, we get that $\chi_\epsilon$ converges uniformly to $\chi$ as $\epsilon \to 0$ on compact subsets in $\R$.   The desired assertion $(i)$ is  \cite[Lemma 5.6]{Viet-generalized-nonpluri}. We consider now $\chi \in \mathcal{W}_M^+$.   By definition, 
$$\chi_\epsilon(t)= \int_\R \chi(t-s) g_\epsilon(s) ds= \int_{\R^+}\big( \chi(t-s)+ \chi(t+s)\big) g_\epsilon(s) ds.$$ 
Hence $\chi_\epsilon$ is concave and increasing, and $\chi_\epsilon \le \chi$.  We note that $\chi'(t)$ is well-defined almost everywhere and decreasing. Hence, 
$$\chi'_\epsilon(t)= \int_\R \chi'(t-s) g_\epsilon(s) ds.$$
Let $t\le 0$ such that $|t|> \epsilon$. Using the fact $\chi \in \mathcal{W}_M^+$ and  that the support of $g_\epsilon$ is contained in the disk of radius $\epsilon$ centered at $0$ gives
\begin{align*}
|t \chi'_\epsilon(t)|  &\le \int_\R  |t\chi'(t-s)| g_\epsilon(s) ds \\
&\le  \int_\R  |(t-s)\chi'(t-s)| g_\epsilon(s) ds+  \int_\R  |s| |\chi'(t-s)| g_\epsilon(s) ds \\
&\le  M \int_\R |\chi(t-s)| g_\epsilon(s) ds+ \epsilon (|t|- \epsilon)^{-1} \int_\R  |(t-s)\chi'(t-s)| g_\epsilon(s) ds \\
&\le -(M+ \epsilon/ (|t|- \epsilon))  \int_\R \chi(t-s) g_\epsilon(s) ds =  -(M+ \epsilon/ (|t|- \epsilon)) \chi_\epsilon(t).
\end{align*}
This finishes the proof.
\endproof

\begin{lemma} \label{le-tau} For every constant $c\ge 1$ and every $t \le 0$, we have 
$$\chi(c t) \ge  c \chi(t)$$
for $\chi \in \mathcal{W}^-$, and
$$\chi(c t) \ge c^M \chi(t)$$
for $\chi \in \mathcal{W}_M^+$.
\end{lemma}
\proof We first assume that $\chi$ is smooth. Consider the function $f(t):= \chi(c t)- c \chi(t)$. We have $f'(t)= c \chi'(c t) - c \chi'(t) \le 0$ because $c\ge 1$, $t\le 0$ and $\chi'$ is increasing. Thus, $f$ is decreasing. In particular, $f(t) \ge f(0)=0$ for $t\le 0$. The first desired inequality hence follows. To prove the second one, since $\chi \in \mathcal{W}_M^+$ and $\chi(t) \le  0$ for $t \le 0$, we get 
$$\frac{\chi'(t)}{\chi(t)} \ge \frac{M}{t} \cdot$$
Integrating both sides along the interval $[ct,t]$ gives the second desired inequality.

Consider now the general case. Let $\chi_\epsilon$ be as in Lemma \ref{le-regularizedchi}. When $\chi \in \mathcal{W}^-$, applying the previous case for $\chi_\epsilon - \chi_\epsilon(0)$ instead of $\chi$ and taking $\epsilon \to 0$ gives the desired inequality for $\chi \in \mathcal{W}^-$. When $\chi \in \mathcal{W}^+_M$, it is essentially the same: we  apply the proof of  the smooth case to $\chi_\epsilon$ and $t<-\delta$ for some fixed constant $\delta>0$, then we get the desired inequality by letting $\epsilon \to 0$ and then $\delta \to 0$. This finishes the proof. 
\endproof

Let $X$ be a compact K\"ahler manifold of dimension $n$ and let $\omega$ be a K\"ahler form on $X$. 
We recall the following result which is taken from the proofs of \cite[Lemmas 2.3 and 3.5]{GZ-weighted}. 

\begin{lemma} \label{le-monotonicityDchi} Let $\varphi_1$ and $\varphi_2$ be bounded negative $\omega$-psh functions on $X$ such that $\varphi_1 \le \varphi_2$. Let $T$ be a closed positive current of bi-dimension $(1,1)$. Then we have
\begin{align}\label{ine-monoreferedoi}
-\int_X \chi(\varphi_1) (\ddc \varphi_2+ \omega) \wedge T \le  - (M+1)\int_X \chi(\varphi_1) (\ddc \varphi_1+ \omega) \wedge T
\end{align}
if $\chi \in \mathcal{W}_M^+$, and 
\begin{align}\label{ine-monoreferedoi2}
-\int_X \chi(\varphi_1) (\ddc \varphi_2+ \omega) \wedge T \le  - 2\int_X \chi(\varphi_1) (\ddc \varphi_1+ \omega) \wedge T
\end{align}
if $\chi \in \mathcal{W}^-$.  
\end{lemma}

\proof
 We reproduce detailed arguments for readers' convenience. By Lemma \ref{le-regularizedchi}, without loss of generality, we can assume that $\chi$ is smooth. Since $\chi$ is increasing, one has $\chi' \ge 0$. Put $\omega_{\varphi_j}:= \ddc \varphi_j+ \omega$ for $j=1,2$. By integration by parts, we have
\begin{align}\label{ine-refereeechiomega} 
- \int_X \chi(\varphi_1) \omega_{\varphi_1} \wedge T= \int_X \chi'(\varphi_1) d \varphi_1 \wedge \dc \varphi_1 \wedge T -\int_X \chi(\varphi_1)\omega\wedge T \ge -\int_X \chi(\varphi_1)\omega\wedge T
\end{align}
(see  \cite[Lemma 5.7]{Viet-generalized-nonpluri} for a comment on the last equality).  Using integration by parts again gives
\begin{align}\label{eq-inrephi12refere}
-\int_X \chi(\varphi_1) \omega_{\varphi_2} \wedge T= -\int_X \varphi_2 \ddc \chi(\varphi_1) \wedge T-\int_X \chi(\varphi_1)  \omega  \wedge T.
\end{align}
Consider now $\chi \in \mathcal{W}^-$. Since $\chi'' \ge 0$, observe that  $\chi'(\varphi_1)\omega+ \ddc \chi(\varphi_1) \ge 0$. By this, the hypothesis and (\ref{eq-inrephi12refere}), we get
\begin{align*}
-\int_X \chi(\varphi_1) \omega_{\varphi_2} \wedge T & \le   -\int_X \varphi_2\big(\chi'(\varphi_1)\omega+  \ddc \chi(\varphi_1)\big) \wedge T- \int_X \chi(\varphi_1)  \omega  \wedge T\\
& \le   -\int_X \varphi_1\big(\chi'(\varphi_1)\omega+  \ddc \chi(\varphi_1)\big) \wedge T- \int_X \chi(\varphi_1)  \omega  \wedge T\\
&=  -\int_X \varphi_1 \chi'(\varphi_1)\omega\wedge T -\int_X \chi(\varphi_1)  \omega_{\varphi_1} \wedge T\\
&\le   -\int_X \chi(\varphi_1)\omega\wedge T -\int_X \chi(\varphi_1)  \omega_{\varphi_1} \wedge T
 \end{align*}
 because $- \varphi_1 \chi'(\varphi_1) \le - \chi(\varphi_1)$ (by convexity of $\chi$ and the equality $\chi(0)=0$). Combining this with (\ref{ine-refereeechiomega}) gives (\ref{ine-monoreferedoi2}). 

It remains to check (\ref{ine-monoreferedoi}) for $\chi \in \mathcal{W}^+_M$. Since $\chi'' \le 0 \le \chi'$, we obtain
$$\ddc \chi(\varphi_1)= \chi''(\varphi_1) d \varphi_1 \wedge \dc \varphi_1+ \chi'(\varphi_1) \ddc \varphi_1 \le \chi'(\varphi_1) \omega_{\varphi_1}.$$
Thus, 
$$-\int_X \varphi_2 \ddc \chi(\varphi_1) \wedge T \le -\int_X \varphi_2  \chi'(\varphi_1) \omega_{\varphi_1}\wedge T \le  -\int_X \varphi_1  \chi'(\varphi_1) \omega_{\varphi_1}\wedge T$$
which is less than or equal to 
$$-M \int_X \chi(\varphi_1) \omega_{\varphi_1} \wedge T$$
because $\chi \in \mathcal{W}^+_M$. Combining this with (\ref{ine-refereeechiomega}) and (\ref{eq-inrephi12refere}) gives  (\ref{ine-monoreferedoi}). The proof is finished.
\endproof

Using Lemma \ref{le-monotonicityDchi} and \cite[Proposition 1.4]{GZ-weighted}, we obtain the following known fact.

\begin{proposition}\label{pro-monoeneryclass} Let $\chi \in \mathcal{W}^- \cup \mathcal{W}^+_M$. Let $\varphi_1$ and $\varphi_2$ be bounded negative $\omega$-psh functions on $X$ such that $\varphi_1 \le \varphi_2$. Then $\varphi_2 \in \mathcal{E}_\chi$ if $\varphi_1$ is so. Furthermore, we have
$$E_\chi(\varphi_2) \le (M+1)^n E_\chi(\varphi_1)$$
for $\chi \in \mathcal{W}^+_M$, and 
$$E_\chi(\varphi_2) \le  2^n E_\chi(\varphi_1)$$
for $\chi \in \mathcal{W}^-$.
\end{proposition}

The following key result  will be proved using an idea from \cite{Viet-convexity-weightedclass,Viet-generalized-nonpluri}. This is the idea used to prove the convexity of weighted classes of currents of relative full mass intersection there. In this section, if no confusion arises,  we will remove the bracket $\langle \quad \rangle$ from the notation of non-pluripolar products for simplicity. 

\begin{theorem} \label{th-integra} There exists a constant $C$ depending only on the dimension $n$ such that for every $\psi,\varphi\in \mathcal{E}_\chi$ such that $\sup_X \psi = \sup_X \varphi =-1$, we have 
\begin{align}\label{eq-intebyparts}
\int_X - \chi(\psi) (\ddc \varphi+ \omega)^n \le  C\bigg(\big(E_{\chi}(\varphi) E_\chi(\psi)\big)^{1/2} +E_{\chi}(\varphi) + \|\chi(\psi)\|_{L^1}\bigg) 
\end{align}
if $\chi \in \mathcal{W}^-$, and 
\begin{align}\label{eq-intebypartsWM}
\int_X - \chi(\psi) (\ddc \varphi+ \omega)^n \le  C (M+1)^n \bigg((E_{\chi}(\varphi))^{1/(M+1)} (E_\chi(\psi))^{M/(M+1)}+E_{\chi}(\varphi) +\|\chi(\psi)\|_{L^1}\bigg)
\end{align}
if $\chi \in \mathcal{W}_M^+$. 
\end{theorem}

We refer to \cite[Lemmas 2.11 and 3.9]{GZ-weighted} for related results  in the case where $\chi(t)= -(-t)^p$ ($p>0$). If $\chi(t)= -(-t)^p$ with $0<p\le 1$, then the exponent $1/2$ in (\ref{eq-intebyparts}) can be replaced by $p/(p+1)$ as the proof of  (\ref{eq-intebypartsWM}) provided below shows. The normalization $\sup_X \psi =-1$ is not essential but it has an advantage that $E_\chi(\psi) \ge \chi(-1)>0$, similarly for $\varphi$.

\proof  We first prove (\ref{eq-intebyparts}).  Let $c\ge 2$ be a constant. Observe 
$$X = \{\psi > c \varphi\} \cup \{\psi < c \varphi+1\}.$$
Hence, 
\begin{align*}
\int_X - \chi(\psi) (\ddc  \varphi+ \omega)^{n} \le \int_{\{\psi> c \varphi\}} - \chi(\psi) (\ddc  \varphi+ \omega)^n +\int_{\{\psi< c \varphi+1\}} - \chi(\psi) (\ddc  \varphi+ \omega)^n.
\end{align*}
Denote by $I_1$ and  $I_2$ the first and second term in the right-hand side of the last inequality. By Lemma \ref{le-tau}, we get 
\begin{align}\label{ine-chanI1Wtru}
I_1\le  \int_{\{\psi> c \varphi\}} - \chi(c \varphi) (\ddc  \varphi+ \omega)^n \le \int_X - \chi(c \varphi) (\ddc  \varphi+ \omega)^n \le c E_\chi(\varphi).
\end{align}
On the other hand, for $\psi':= \max\{\psi, c \varphi +1\} \ge \psi$ (which is a $(c\omega)$-psh function), using the plurifine locality, one has
\begin{align*}
I_2  &= c^{-n}\int_{\{\psi< c \varphi+1\}} - \chi(\psi) (\ddc  (c\varphi)+ c \omega)^{n}\\
& = c^{-n}\int_{\{\psi< c \varphi+1\}} - \chi(\psi) (\ddc  \psi'+ c \omega)^{n} \le c^{-n}\int_X - \chi(\psi) (\ddc  \psi'+ c \omega)^{n} 
\end{align*}
which is, by monotonicity of energy classes (see Lemma \ref{le-monotonicityDchi}),  less than or equal to
$$2^n c^{-n} \int_{X} - \chi(\psi) (\ddc  \psi+ c \omega)^{n}.$$
The last term is less than or  equal to
$$ A c^{-n} (c-1)^{n} \|\chi(\psi)\|_{L^1}+ A c^{-n} \sum_{j=1}^n \int_{X} - \chi(\psi) (\ddc  \psi+ \omega)^{j} \wedge [(c-1)\omega]^{n-j},$$
where $A\ge 1$ is a dimensional constant. The second term in the last sum is bounded by 
\begin{align}\label{ine-ActrungWtru}
A c^{-n} \sum_{j=1}^n c^{n-j} \int_X  - \chi(\psi) (\ddc  \psi+ \omega)^{j} \wedge \omega^{n-j} \le   n A c^{-1} E_{\chi}(\psi)
\end{align}
 by Lemma \ref{le-monotonicityDchi} again. It follows that there exists a constant $A'>0$ depending only on $n$ such that
\begin{align}\label{ine-uocluongmixedenergy}
\int_X - \chi(\psi) (\ddc  \varphi+ \omega)^{n}  \le A' \big(c E_\chi(\varphi)+ c^{-1}E_\chi(\psi)+\|\chi(\psi)\|_{L^1}\big)
\end{align}
for every constant  $c \ge 2$. Let $c_0:= (E_\chi(\psi)/ E_\chi(\varphi))^{1/2}$. If $c_0\le 2$, then  (\ref{eq-intebyparts}) is clear. Consider $c_0 \ge 2$. Applying (\ref{ine-uocluongmixedenergy}) to $c_0$ yields the desired assertion. 

The other inequality (\ref{eq-intebypartsWM}) is proved similarly. We only indicate minor changes needed to implement. Consider now $\chi \in\mathcal{W}_M^+$.  In the right-hand sides of (\ref{ine-chanI1Wtru}) and (\ref{ine-ActrungWtru}), the factors $c$ and $ n A c^{-1}$  are replaced by $c^M$ and $(M+1)^n n A c^{-1}$, respectively. As a result, we get
$$\int_X - \chi(\psi) (\ddc  \varphi+ \omega)^{n}  \le A'(M+1)^n  \big(c^M E_\chi(\varphi)+ c^{-1}E_\chi(\psi)+\|\chi(\psi)\|_{L^1}\big).$$
Arguing as before gives   (\ref{eq-intebypartsWM}).   This finishes the proof.
\endproof

Note that the term $\|\chi(\psi)\|_{L^1}$ is bounded by a constant independent of $\psi$ with $\sup_X \psi =-1$ because $|\chi(\psi)| \le - \chi(-1)|\psi| $ if $\chi \in \mathcal{W}^-$, and $|\chi(\psi)| \le -\chi(-1)|\psi|^M $ if $\chi \in \mathcal{W}_M^+$ (see Lemma \ref{le-tau}).
The following is a well-known result from \cite{GZ-weighted}, see also \cite{Cegrell}.

\begin{theorem}\label{th-dominatedcapacity} Let $\mu$ be a probability measure dominated by capacity, \emph{i.e,} there exists a constant $A$ such that for every Borel set $K$ in $X$, we have $\mu(K) \le A \, \capK_\omega(K)$. Then there exists an $\omega$-psh function $\varphi \in \cap_{p \ge 1} \mathcal{E}^p$ such that $\mu= (\ddc \varphi+ \omega)^n$.
\end{theorem}

\proof  One can see Theorem \ref{th-dominatedcapacitybig} below for a more general statement for big cohomology classes. Although the desired result was not explicitly stated in \cite{GZ-weighted}, it was indeed implicitly proved in the proof of Theorem  4.1 in \cite{GZ-weighted}.  Firstly, as one can see from the paragraph right before \cite[Lemma 4.5]{{GZ-weighted}}, there exists $\varphi \in \mathcal{E}^1$ such that $\mu= (\ddc \varphi+ \omega)^n$.  Now by the last paragraph of the proof of \cite[Theorem 4.1]{GZ-weighted} (the end of the page 472 and the beginning of the page 473 there), one gets $\varphi \in \cap_{p \ge 1} \mathcal{E}^p$. 
\endproof

Recall that 
$$\capK_\omega(K):= \sup \big\{ \int_K (\ddc \psi+ \omega)^n: \, \psi \quad \text{$\omega$-psh}, \quad  0 \le \psi \le 1\big\}.$$
Note that since $|\chi(t)| \le -|t|^M \chi(-1)$ for $t \le -1$, we have 
$$\cap_{p \ge 1} \mathcal{E}^p \subset \mathcal{E}_\chi$$
for every $\chi \in \mathcal{W}_M^+$. It was shown in \cite{Kolodziej_Acta} that the solution $\varphi$ in Theorem \ref{th-dominatedcapacity} is not bounded in general.


\begin{proof}[End of the proof of Theorem \ref{th-main}]  Firstly, assume that $\mu= (\ddc \varphi+ \omega)^n$ for some $\varphi \in \mathcal{E}_\chi$. Applying Theorem \ref{th-integra} to $\mu$ gives the desired inequality. 

Now suppose (\ref{ine-dktuongduongEchi}) is satisfied. We will prove the desired assertion under a weaker hypothesis that there exists a function $F: \R^+ \to \R^+$  with $\limsup_{t \to \infty} F(t)/t <1$ and 
$$-\int_X \chi(\psi) d \mu \le F\big(E_\chi(\psi)\big)$$
for every  smooth $\omega$-psh function $\psi$ with $\sup_X \psi=-1$.  It is clear that $\mu$ has no mass on pluripolar sets. Thus, by \cite[Theorem 4.6]{GZ-weighted}, there exists $\varphi\in \mathcal{E}$ so that $\mu= (\ddc \varphi+ \omega)^n$ (this solution is unique by \cite{Dinew-uniqueness}). We need to check that $\varphi \in \mathcal{E}_\chi$. However we don't really need \cite[Theorem 4.6]{GZ-weighted} here because the existence of $\varphi$ will be deduced from the arguments in the next paragraph.  The strategy is  that of  \cite{Cegrell} and \cite[Theorem 4.1]{GZ-weighted}. A key point in this approach is to control the energy of the sequence $(\varphi_j)_j$ below. 

We recall how to construct $\varphi$. Since $\mu$ has no mass on pluripolar sets,  there exists a probability measure $\nu$ bounded by capacity and a Borel function $f\ge 0$ locally integrable with respect to $\nu$ such that $\mu= f \nu$; see for example \cite[Lemma 4.5]{GZ-weighted}.  Let 
$$\mu_j:= \delta_j \min \{f, j\} \nu,$$
where $\delta_j>0$ is a constant such that the mass of $\mu_j$ is equal to $1$. Note that $\mu_j$ converges to $\mu$ as $j \to \infty$, and $\mu_j \le \delta_j \mu$. Observe also that $\delta_j$ converges to $1$ as $j \to \infty$.  By Theorem \ref{th-dominatedcapacity}, there exists $\varphi_j \in \cap_{p \ge 1} \mathcal{E}^p$ satisfying 
$$\mu_j= (\ddc \varphi_j + \omega)^n, \quad \sup_X \varphi_j =-1.$$
By (\ref{ine-dktuongduongEchi}), we have
$$E_\chi(\varphi_j)= - \int_X \chi(\varphi_j) (\ddc \varphi_j+ \omega)^n \le  -\delta_j \int_X \chi(\varphi_j) d \mu \le \delta_j F\big(E_\chi(\varphi_j)\big).$$
This combined with the fact that $\limsup_{t \to \infty}F(t)/t <1$ gives  
\begin{align}\label{ine-supenergyvapghij}
\sup_j E_\chi(\varphi_j) < \infty.
\end{align} 
By extracting a subsequence if necessary, we can assume that $\varphi_j \to \varphi$ in $L^1$ for some $\omega$-psh function $\varphi$. Hence, $\varphi$ is the limit of the decreasing sequence of $\varphi_j':= (\sup_{k \ge j} \varphi_j)^*$. Note that since $\varphi'_j \ge \varphi_j$, we get $\varphi'_j \in \mathcal{E}_\chi$ (monotonicity of energy classes, Proposition \ref{pro-monoeneryclass}) and 
$$\sup_j E_{\chi}(\varphi'_j) \lesssim \sup_j E_\chi(\varphi_j)+1 <\infty$$ 
 by  (\ref{ine-supenergyvapghij}). By the lower semi-continuity like property of energy (see \cite[Corollary 2.7]{GZ-weighted} or Corollary \ref{cor-semi-continuity-big} below for general big cohomology case), we obtain that 
 $$E_\chi(\varphi) \lesssim \liminf_{j \to \infty} E_\chi(\varphi'_j) < \infty$$
 and $\varphi$ is of full Monge-Amp\`ere mass. In other words, $\varphi \in \mathcal{E}_\chi$. It follows that 
$$(\ddc \varphi + \omega)^n = \lim_{j \to \infty} (\ddc \varphi'_j + \omega)^n \ge \limsup_{j \to \infty} \inf_{k \ge j} \min \{f,k\} \nu \ge f \nu= \mu.$$
Since both sides of the above inequality has the same mass on $X$, we get 
$$(\ddc \varphi+ \omega)^n =\mu.$$
So $\varphi \in \mathcal{E}_\chi$ is what we are looking for. The proof is finished.
\end{proof}

\section{Big cohomology class setting} \label{sec-bigcoho}

In this section, we explain how to extend Theorem \ref{th-main} to the case of big cohomology class.  Let $\theta$ be a closed smooth  real $(1,1)$-form on $X$.  Let $\PSH(X, \theta)$ be the set of $\theta$-psh functions on $X$.  Let $\phi$ and $\psi$ be in $\PSH(X, \theta)$. Recall that $\psi$ is said to be \emph{more singular than} $\phi$ if $\psi \le \phi + O(1)$, in this case, we write $\psi \preceq \phi$. The function $\phi$ and $\psi$ are said to be of \emph{the same singularity type}  if $\phi$ is more singular than $\psi$ and vice versa.

 We first extend Lemma \ref{le-monotonicityDchi} to our setting. This was essentially done in \cite{Viet-convexity-weightedclass}. We provide details for readers' convenience. 
To do so, we need an integration by parts formula for relative non-pluripolar products from \cite{Viet-convexity-weightedclass} generalizing those given in \cite{BEGZ,Lu-cap-compare,Xia}. We will employ the notion of relative non-pluripolar products in \cite{Viet-generalized-nonpluri}. This use is not absolutely necessary in the proof of Theorem \ref{th-mainbigclass} but it makes the presentation more clear.

Let $T_1, \ldots, T_m$ be closed positive $(1,1)$-currents on $X$. Let $T$ be  a closed positive current of bi-degree $(p,p)$ on $X$. The \emph{$T$-relative non-pluripolar product} $\langle \wedge_{j=1}^m T_j \dot{\wedge} T\rangle$ is defined  in a way similar to that of  the usual non-pluripolar product. The product $ \langle  \wedge_{j=1}^m T_j \dot{\wedge} T\rangle $ is a well-defined closed positive current of bi-degree $(m+p,m+p)$; and $ \langle  \wedge_{j=1}^m T_j \dot{\wedge} T\rangle $  is  symmetric with respect to $T_1, \ldots, T_m$ and is homogeneous. In latter applications, we will only use the case where $T$ is the non-pluripolar product of some closed positive $(1,1)$-currents, say, $T= \langle T_{m+1} \wedge \cdots \wedge T_{m+l} \rangle$, where $T_j$ is $(1,1)$-currents for $m+1 \le j \le m+l$. In this case, $\langle T_1 \wedge \cdots \wedge T_m \dot{\wedge}T \rangle $ is simply equal to $\langle \wedge_{j=1}^{m+l} T_j \rangle$. 

Recall that a \emph{dsh} function on $X$ is the difference of two quasi-plurisubharmonic (quasi-psh for short) functions on $X$ (see \cite{DS_tm}). These functions are well-defined outside pluripolar sets. Let $v$ be a dsh function on $X$.  Let $T$ be a closed positive current on $X$. We say that $v$ is \emph{$T$-admissible} if  there exist  quasi-psh functions $\varphi_1, \varphi_2$ such that $v= \varphi_1- \varphi_2$  and $T$ has no mass on $\{\varphi_j=-\infty\}$ 
for $j=1,2$. In particular, if $T$ has no mass on pluripolar sets, then every dsh function is $T$-admissible.  

Assume now that $v$ is  \emph{bounded} $T$-admissible.    Let $\varphi_{1}, \varphi_{2}$ be quasi-psh functions such that $v= \varphi_{1}- \varphi_{2}$ and $T$ has no mass on $\{\varphi_{j}=-\infty\}$ for $j=1,2$. Let 
$$\varphi_{j,k}:= \max\{\varphi_{j}, -k \}$$
for every $j=1,2$ and $k \in \N$. Put $v_k:= \varphi_{1,k}- \varphi_{2,k}$. Put
$$Q_k:= d v_k \wedge \dc v_k \wedge T=\ddc v_k^2  \wedge T - v_k\ddc v_k \wedge T.$$
By the plurifine locality with respect to $T$ (\cite[Theorem 2.9]{Viet-generalized-nonpluri}) applied to the right-hand side of the last equality, we have 
\begin{align}\label{eq-localplurifineddc}
\bold{1}_{\cap_{j=1}^2 \{\varphi_{j}> -k\}} Q_k =\bold{1}_{\cap_{j=1}^2\{\varphi_{j}> -k\}} Q_{k'}
\end{align}
for every $k'\ge k$. By \cite[Lemma 2.5]{Viet-convexity-weightedclass},  the mass of $Q_k$ on $X$ is bounded uniformly in $k$.  This combined with  (\ref{eq-localplurifineddc}) implies that there exists a positive current $Q$ on $X$ such that for every bounded Borel form $\Phi$ with compact support on $X$ such that 
$$\langle Q,\Phi \rangle  = \lim_{k\to \infty} \langle Q_k, \Phi\rangle.$$
We define $\langle d v \wedge \dc v \dot{\wedge} T \rangle$ to be the current $Q$.  This agrees with the classical definition if $v$ is the difference of two  bounded quasi-psh functions. One can check  that this definition is independent of the choice of $\varphi_1, \varphi_2$.

Let $w$ be another bounded $T$-admissible dsh function.  If $T$ is of bi-degree $(n-1,n-1)$, we can also define the current $\langle  dv \wedge \dc w \dot{\wedge} T\rangle$ by a similar procedure as above.  We put
$$\langle \ddc v \dot{\wedge} T \rangle:= \langle \ddc \varphi_1 \dot{\wedge} T\rangle- \langle \ddc \varphi_2 \dot{\wedge} T\rangle$$
which is independent of the choice of $\varphi_1,\varphi_2$.  The following integration by parts formula is crucial for us later.

\begin{theorem} \label{th-integrabypart}  (\cite[Theorem 2.6]{Viet-convexity-weightedclass}) Let  $T$ a closed positive current of bi-degree $(n-1,n-1)$ on $X$. Let $v,w$ be bounded $T$-admissible dsh functions on $X$. Let $\chi: \R \to \R$ be a $\cali{C}^3$ function. Then we have  
\begin{multline}\label{eq-intebypartschi}
\int_X \chi(w) \langle  \ddc v \dot{\wedge} T\rangle=\int_X v \chi''(w) \langle dw \wedge \dc w \dot{\wedge} T\rangle+\int_X v \chi'(w) \langle \ddc w \dot{\wedge} T\rangle\\
= - \int_X \chi'(w) \langle d w \wedge  \dc v \dot{\wedge} T\rangle.
\end{multline}
\end{theorem}

Note that the second equality of (\ref{eq-intebypartschi}) wasn't stated explicitly in \cite[Theorem 2.6]{Viet-convexity-weightedclass}. But its validity can be seen directly from the proof of the last result.  Here is a generalization of Lemma \ref{le-monotonicityDchi} to our setting. We note that some partial results were obtained in \cite[Section 2.2]{Lu-Darvas-DiNezza-logconcave}.

\begin{lemma} \label{le-monotonicityDchibig} Let $\varphi_1,\varphi_2$ and $\phi$ be $\theta$-psh functions on $X$ such that $\varphi_1,\varphi_2$ are more singular than $\phi$, and $\varphi_1 \le \varphi_2$. Let $T$ be a closed positive current of bi-dimension $(1,1)$ having no mass on pluripolar sets. Assume that 
\begin{align}\label{eq-varphi1fullmassT}
\int_X \langle (\ddc \varphi_1+ \theta)  \dot{\wedge} T\rangle= \int_X \langle (\ddc \phi+ \theta)  \dot{\wedge} T\rangle.
\end{align}
Then we have
$$-\int_X \chi(\varphi_1-\phi) \langle (\ddc \varphi_2+ \theta) \dot{\wedge} T\rangle \le  - (M+1)\int_X \chi(\varphi_1-\phi) \langle (\ddc \varphi_1+ \theta)  \dot{\wedge} T\rangle$$ 
if $\chi \in \mathcal{W}_M^+$, and 
$$-\int_X \chi(\varphi_1-\phi) \langle (\ddc \varphi_2+ \theta) \dot{\wedge} T\rangle \le  - 2\int_X \chi(\varphi_1-\phi) \langle (\ddc \varphi_1+ \theta)  \dot{\wedge} T\rangle $$
if $\chi \in \mathcal{W}^-$.  
\end{lemma}

\proof By regularization, we can assume that $\chi$ is smooth.  The case $\chi \in \mathcal{W}^-$ was treated in \cite{Viet-convexity-weightedclass,Viet-generalized-nonpluri}. The case $\chi \in \mathcal{W}_M^+$ is also essentially the same as the last one. We explain briefly how to proceed.  

Using (\ref{eq-varphi1fullmassT}), the fact that $\varphi_1 \le \varphi_2$ and the monotonicity of relative non-pluripolar products (\cite[Theorem 1.1]{Viet-generalized-nonpluri}) gives that 
\begin{align}\label{eq-varphi1fullmassT2}
\int_X \langle (\ddc \varphi_2+ \theta)  \dot{\wedge} T\rangle= \int_X\langle (\ddc \phi+ \theta)  \dot{\wedge} T\rangle.
\end{align}
Put $\varphi_{j,k}:= \max\{\varphi_j- \phi, -k\}$ for $j=1,2$. We see that $\varphi_{j,k}$ is a bounded $T$-admissible dsh function converging to $\varphi_j- \phi$ as $k \to \infty$.  

By \cite[Lemma 5.2]{Viet-generalized-nonpluri}, (\ref{eq-varphi1fullmassT}) and (\ref{eq-varphi1fullmassT2}), in order to prove the desired inequality, it suffices to prove it with $\varphi_{j,k}$ in place of $\varphi_j- \phi$ (hence $\chi(\varphi_j- \phi)$ and $\ddc \varphi_j+ \theta$  are substituted by $\chi(\varphi_{j,k})$ and $\ddc \varphi_{j,k}+ \theta_\phi$ respectively, where $\theta_\phi:= \ddc \phi+ \theta$).  

Now, one just needs to follow literally arguments in the proof of Lemma \ref{le-monotonicityDchi} presented in \cite{GZ-weighted} (with $\omega$ replaced by $\theta_\phi$ and the classical product replaced by the non-pluripolar one) and use the integration by parts formula given by Theorem \ref{th-integrabypart} (we note however that  other weaker integration by parts formulae obtained in \cite{BEGZ,Lu-cap-compare,Xia} seem to be not enough  for our purposes here). This finishes the proof. 
\endproof

Let $\mathcal{E}(\theta, \phi)$ be the set of $\varphi \in \PSH(X, \theta)$ such that $\varphi$ is more singular than $\phi$ and 
$$\int_X \big \langle (\ddc \varphi+\theta)^n \big \rangle = \int_X \big \langle (\ddc \phi+\theta)^n \big \rangle.$$
The last space was introduced  in \cite{Lu-Darvas-DiNezza-mono}, see also \cite{Viet-convexity-weightedclass} for some generalizations.   We recall the notion of $\chi$-energy. The case $\chi(t)=t$ was studied in \cite{Lu-Darvas-DiNezza-logconcave}. For $\varphi \in \mathcal{E}(\theta, \phi)$, we put 
$$E_{\chi,\theta, \phi}(\varphi):=- \int_X \chi(\varphi - \phi) \big \langle (\ddc \varphi+ \theta)^n\big \rangle.$$
When $\phi$ is of minimal singularities in the cohomology class of $\theta$, the above energy is equivalent to that defined in \cite{Viet-convexity-weightedclass}. Furthermore, if $\phi' \in \PSH(X,\theta)$ such that $\phi'$ and $\phi$ are of the same singularity type, then by Lemma \ref{le-tau} we have 
\begin{align}\label{ine-Engersingulartype}
C^{-1 } E_{\chi,\theta,\phi'}(\varphi)  \le  E_{\chi,\theta,\phi}(\varphi) \le C E_{\chi,\theta,\phi'}(\varphi)
\end{align}
for some constant $C>0$ independent of $\varphi$.

 Let $\mathcal{E}_{\chi, \theta,\phi}$ be the subset of $\mathcal{E}(\theta, \phi)$ consisting of $\varphi$ with $E_{\chi, \theta,\phi}(\varphi)<\infty$. When $\chi(t)= -(-t)^p$, we write $E_{p,\theta,\phi}$, $\mathcal{E}^p_{\theta,\phi}$ for  $E_{\chi,\theta,\phi}$, $\mathcal{E}_{\chi,\theta,\phi}$ respectively.   Using Lemma \ref{le-monotonicityDchibig}, we get the following important estimates just as in the K\"ahler case.

\begin{proposition}\label{pro-monoeneryclassbig}Let $\chi \in \mathcal{W}^- \cup \mathcal{W}^+_M$. Let $\varphi_1$ and $\varphi_2$ be in $\PSH(X, \theta)$ such that $\varphi_1 \le \varphi_2$. Assume that  $\varphi_1 \in \mathcal{E}_{\chi, \theta, \phi}$. Then $\varphi_2 \in \mathcal{E}_{\chi, \theta, \phi}$, and we have
$$E_{\chi,\theta,\phi}(\varphi_2) \le (M+1)^n E_{\chi,\theta,\phi}(\varphi_1)$$
for $\chi \in \mathcal{W}^+_M$, and 
$$E_{\chi,\theta,\phi}(\varphi_2) \le  2^n E_{\chi,\theta,\phi}(\varphi_1)$$
for $\chi \in \mathcal{W}^-$.
\end{proposition}

As a consequence, we obtain a lower semi-continuity like property of energy. 

\begin{corollary} \label{cor-semi-continuity-big} Let $\chi \in \mathcal{W}^- \cup \mathcal{W}^+_M$.  Let $(\varphi_k)_k$ be a sequence of $\theta$-psh functions in $\mathcal{E}_{\chi,\theta,\phi}$ converging to a $\theta$-psh function $\varphi$ in $L^1$ such that $E_{\chi,\theta,\phi}(\varphi_k)$ is bounded uniformly in $k$. Then we have $\varphi \in \mathcal{E}_{\chi,\theta,\phi}$, and  
$$E_{\chi,\theta,\phi}(\varphi) \le 2^n \liminf_{k \to \infty} E_{\chi,\theta,\phi}(\varphi_k)$$
 if $\chi \in \mathcal{W}^-$, and 
$$E_{\chi,\theta,\phi}(\varphi) \le (M+1)^n \liminf_{k \to \infty} E_{\chi,\theta,\phi}(\varphi_k)$$
 if $\chi \in \mathcal{W}^+_M$.
\end{corollary}

\proof We consider $\chi \in \mathcal{W}^-$. The other case is treated similarly. Let $\varphi'_k:= (\sup_{l\ge k} \varphi_l)^*$. We have that $\varphi'_k$ decreases pointwise to $\varphi$ as $k \to \infty$. Using this and the fact that $\chi$ is continuous increasing gives
$$E_{\chi,\theta,\phi}(\varphi) \le \liminf_{k\to \infty} E_{\chi,\theta,\phi}(\varphi'_k).$$
 By Proposition \ref{pro-monoeneryclassbig}, we get 
$$E_{\chi,\theta,\phi}(\varphi'_k) \le 2^n \inf_{l\ge k} E_{\chi,\theta,\phi}(\varphi_l)$$
for every $k$. Letting $k\to \infty$ gives $E_{\chi,\theta,\phi}(\varphi)<\infty$. It remains to check that $\varphi \in \mathcal{E}(\theta,\phi)$. Observe that by above arguments applied to $\psi_l:= \max\{\varphi, \phi -l\}$ for $l\ge 1$ gives 
$$E_{\chi,\theta,\phi}(\psi_l) \lesssim \liminf_{k \to\infty} E_{\chi,\theta,\phi}(\max\{\varphi_k,\phi -l\}) \lesssim \liminf_{k \to\infty} E_{\chi,\theta,\phi}(\varphi_k)$$
by monotonicity. It follows that 
$$\int_{\{\varphi\le  \phi-l\}} \langle (\ddc \psi_l+ \theta)^n \rangle = - \big|\chi(-l)\big|^{-1} \int_{\{\varphi\le  \phi-l\}} \chi(\psi_l -\phi)  \langle (\ddc \psi_l+ \theta)^n \rangle \lesssim \big|\chi(-l)\big|^{-1}$$
which tends to $0$ as $l\to \infty$ because $\chi(-\infty)= -\infty$ and $\chi$ is continuous. We deduce that $\varphi \in \mathcal{E}(\theta,\phi)$.  The proof is finished.
\endproof

We now go back to Monge-Amp\`ere equations. Let $\phi$ and  $\psi$ be negative $\theta$-psh functions on $X$. Recall 
$$P_\theta(\phi, \psi):= \big(\sup \{ v\in \PSH(X, \theta): v \le \min\{\phi, \psi\}\}\big)^*$$
if the set in the supremum is non-empty, and $P_\theta(\phi, \psi):= -\infty$ otherwise. Put
$$P_\theta[\phi]:= \big( \sup\{v \in \PSH(X, \theta): v \le 0,\quad v \preceq \phi\} \big)^*.$$
Observe that $P_\theta(\phi+C, 0)$ increases a.e. to $P_\theta[\phi]$ as $C\in \R$ tends to $\infty$. We always have $\phi \preceq P_\theta[\phi]$. Moreover, the singularities of $\phi$ and $P_\theta[\phi]$ are similar in the sense that the multiplier ideal sheafs associated to them are the same, see \cite[Page 405]{Lu-Darvas-DiNezza-singularitytype}. 

Let $\alpha$ be the cohomology class of $\theta$. A \emph{singularity type} in $\alpha$ is an equivalence class of $\theta$-psh functions of the same singularity type.  Given a singularity type $\xi$ in $\alpha$, we define the volume of $\xi$ by putting
$$\vol(\xi):= \int_X \big \langle (\ddc \phi+\theta)^n \big \rangle,$$  
where $\phi \in \xi$. The last quantity is independent of the choice of $\phi$ by monotonicity.   The singularity type of a $\theta$-psh function is the singularity type in $\alpha$ containing that function.   We say that $\phi$ is of \emph{model type singularity} if $\phi$ and $P_\theta[\phi]$ are of the same singularity type. A \emph{model singularity type} is the singularity type of a $\theta$-psh function of model type singularity.

It was known that $\theta$-psh functions with minimal singularities or with analytic singularity are of model type singularities, see \cite[Proposition 4.36]{Lu-Darvas-DiNezza-mono} or \cite{Rashkovskii-Sigurdsson,Ross-WittNystrom}. By \cite[Theorem 3.12]{Lu-Darvas-DiNezza-mono}, the function $P_\theta[\phi]$ is of model type singularity for every $\phi$ with $\int_X \theta_\phi^n >0$. The model singularity types appear naturally  in the study of complex Monge-Amp\`ere equations. The use of this notion is in fact necessary as pointed out in \cite[Section 4E]{Lu-Darvas-DiNezza-mono}.

Assume now that the cohomology class $\alpha$ of $\theta$ is big. Let $\xi$ be a model singularity type in $\alpha$ with $\vol(\xi)>0$. Let $\phi \in \PSH(X, \theta)$ be in $\xi$.  Using above results and following arguments in the proof  Theorem \ref{th-integra} (to be precise, compared to the proof of the aforementioned result, we need to substitute $\omega$ by $\ddc \phi+\theta $ and $\psi, \varphi$ by $\psi-\phi, \varphi -\phi$ respectively), we obtain the following generalization of  Theorem \ref{th-integra} to the big cohomology class setting.

\begin{theorem} \label{th-integrabig} Let $\mu_\phi:= \big \langle (\ddc \phi+\theta)^n\big \rangle $. There exists a constant $C$ depending only on the dimension $n$ such that for every $\psi,\varphi\in \mathcal{E}_\chi(\theta, \phi)$ such that $\sup_X (\psi-\phi) = \sup_X (\varphi-\phi) =-1$, we have 
\begin{align}\label{eq-intebypartsbig}
\int_X - \chi(\psi- \phi) \big\langle (\ddc \varphi+ \theta)^n \big \rangle  \le  C\bigg(\big( E_{\chi,\theta,\phi}(\varphi)  E_{\chi,\theta,\phi}(\psi)\big)^{1/2} + E_{\chi,\theta,\phi}(\varphi)+  \|\chi(\psi-\phi)\|_{L^1(\mu_\phi)}\bigg) 
\end{align}
if $\chi \in \mathcal{W}^-$, and 
\begin{multline}\label{eq-intebypartsWMbig}
\int_X - \chi(\psi-\phi) \big\langle (\ddc \varphi+ \theta)^n \big \rangle \le \\
 C (M+1)^n \bigg( E_{\chi,\theta,\phi}(\varphi))^{1/(M+1)} ( E_{\chi,\theta,\phi}(\psi))^{M/(M+1)}+ E_{\chi,\theta,\phi}(\varphi)+\|\chi(\psi-\phi)\|_{L^1(\mu_\phi)}\bigg)
\end{multline}
if $\chi \in \mathcal{W}_M^+$. 
\end{theorem}

As a consequence, we get the following result. 

\begin{corollary}\label{cor-integgralbigchonphi} Let $\varphi\in \mathcal{E}_\chi(\theta, \phi)$ with $\sup_X (\varphi-\phi) =-1$. Then there exists a constant $C>0$ such that for every $\psi\in \mathcal{E}_\chi(\theta, \phi)$ such that $\sup_X (\psi-\phi) =-1$, we have 
\begin{align*} 
\int_X - \chi(\psi- \phi) \big\langle (\ddc \varphi+ \theta)^n \big \rangle \le  C\big(\big(E_{\chi,\theta,\phi}(\psi)\big)^{1/2} + 1\big) 
\end{align*}
if $\chi \in \mathcal{W}^-$, and 
\begin{align*}
\int_X - \chi(\psi-\phi) \big\langle (\ddc \varphi+ \theta)^n \big \rangle \le
 C\big(( E_{\chi,\theta,\phi}(\psi))^{M/(M+1)}+ 1\big)
\end{align*}
if $\chi \in \mathcal{W}_M^+$ ($C$ might depend on $M$). 
\end{corollary}

\proof The desired assertions follow from Theorem \ref{th-integrabig} if we can show that $\|\chi(\psi-\phi)\|_{L^1(\mu_\phi)}$ is finite. It is not clear to us if this is true for an arbitrary $\phi$. We overcome this difficulty as follows.  Let  $\phi' \in \PSH(X, \theta)$ be of the same singularity type as $\phi$. We will make clear the choice of $\phi'$ later.   By (\ref{ine-Engersingulartype}) and Lemma \ref{le-tau}, it suffices to prove the desired inequalities for $\phi'$ in place of $\phi$.     

Since $\xi$ is a model singularity type of positive volume, using  \cite[Theorem A]{Lu-Darvas-DiNezza-logconcave}, we obtain that  there is a $\theta$-psh function $\phi'$ in $\xi$ such that $\mu_{\phi'}:= \big \langle (\ddc \phi'+ \theta)^n \big \rangle$ is a smooth measure. 
Hence,  there exist constants $c,C>0$ such that for every $\omega$-psh function $\psi$ with $\sup_X \psi =0$, we have
\begin{align}\label{ine-Skodamodeltypesing}
\int_X e^{-c \psi} \, d \mu_{\phi'} \le C.
\end{align}
By this  and \cite[Lemma 2.2]{Lu-Darvas-DiNezza-logconcave}, the quantity $ \|\chi(\psi-\phi)\|_{L^1(\mu_{\phi'})}$ is bounded uniformly in $\psi$. Hence the desired assertions follow from Theorem \ref{th-integrabig} applied to $\phi'$.
\endproof

We recall that 
$$\capK_{\theta,\phi}(K):= \sup\big\{ \int_K (\ddc \psi+ \theta)^n: \, \psi \in \PSH(X, \theta), \quad 0 \le \psi - \phi \le 1\big\}.$$
The last notion was introduced in  \cite{Lu-Darvas-DiNezza-mono,DiNezzaLu-capacity,DiNezzaLu-quasiprojective} as a generalization of the usual capacity $\capK_\omega$. The following result generalizes Theorem \ref{th-dominatedcapacity}.

\begin{theorem}\label{th-dominatedcapacitybig} Let $\mu$ be a probability measure dominated by capacity, \emph{i.e,} there exists a constant $A$ such that for every Borel set $K$ in $X$, we have $\mu(K) \le A \, \capK_{\theta,\phi}(K)$. Then there exists an $\theta$-psh function $\psi \in \cap_{p \ge 1} \mathcal{E}^p_{\theta,\phi}$ such that $\mu= (\ddc \psi+ \theta)^n$.
\end{theorem}

\proof By \cite[Lemma 2.7]{Lu-Darvas-DiNezza-logconcave} and its proof, for every constant $p>1$, there is a constant $C_p$ independent of $\phi$ such that for every $u \in \mathcal{E}(\theta, \phi)$ with $\sup_X u=0$, there holds 
$$\int_X |u-\phi|^p \, d \mu \le C_p\big(E_{p-1, \theta, \phi}(u)+1\big).$$
Now using  \cite[Theorem 4.7]{Lu-Darvas-DiNezza-logconcave} (or rather its proof), there exists a $\theta$-psh $\psi \in \mathcal{E}^1(\theta, \phi)$ with $\sup_X \psi=0$ such that  $\mu= (\ddc \psi+ \theta)^n$. The above inequality applied to $u:= \psi$ shows that $\psi \in \cap_{p \ge 1} \mathcal{E}^p_{\theta,\phi}$ (by induction). This finishes the proof.
\endproof

\begin{proof}[End of the proof of Theorem \ref{th-mainbigclass}] If $\mu= \big \langle (\ddc \varphi+ \theta)^n\big \rangle$ for some $\varphi \in  \mathcal{E}_{\chi}(\theta, \phi)$, then by Corollary \ref{cor-integgralbigchonphi} we get  (\ref{ine-dktuongduongEchibig}). We prove the converse implication.  With all of above results, the proof is now standard. Firstly, as in the proof of \cite[Theorem 4.7]{Lu-Darvas-DiNezza-logconcave}, there exist a measure $\nu$ bounded by capacity (\emph{i.e,} $\nu(K) \le M\capK_{\theta,\phi}(K)$ for every Borel set $K$ and some constant $M$ independent of $K$) and a nonnegative function $f \in L^1(\nu)$ such that $\mu= f \nu$. We now define $\mu_j$ as in the proof of Theorem \ref{th-main}, and  just follow almost word by word arguments there to get the desired assertion. Here are minor modifications: the form $\omega$ is replaced by $\ddc \phi+\theta$, and  $\varphi_j, \varphi$ are replaced by $\varphi_j -\phi$, $\varphi - \phi$ respectively; Theorem \ref{th-dominatedcapacitybig}, Proposition \ref{pro-monoeneryclassbig} substitute Theorem \ref{th-dominatedcapacity}, Proposition \ref{pro-monoeneryclass}, respectively.    
\end{proof}

As mentioned in Introduction, we now give another proof of the following result using the variational method introduced in \cite{BBGZ-variational}.

\begin{theorem}\label{th-variation} (\cite[Theorem A]{Lu-Darvas-DiNezza-logconcave}) Let $\theta$ and $\alpha$ be as above. Let $\xi$ be a model singularity type in $\alpha$ with $\vol(\xi)>0$. Let $f \in L^p(\omega^n), p > 1$ be
such that $f \ge 0$ and $\int_X f \omega^n = \vol(\xi)$. Then the following hold:

$(i)$ there exists $\varphi \in \xi$, unique up to a constant,  solving (\ref{eq-MA}) for $\mu=f \omega^n$,

$(ii)$ for any constant $\lambda > 0$ there exists a unique $\varphi \in \xi$ such that
$$\big \langle (\ddc \varphi+\theta)^n\big\rangle = e^{\lambda \varphi} f \omega^n.$$
\end{theorem}

\proof
We note that this result was proved in \cite[Section 4]{Lu-Darvas-DiNezza-mono}, using the variational method, under the hypothesis  that $\xi$ has a small unbounded locus. In the general case, we follow literally arguments there. In the list below, we will indicate which results  in \cite[Section 4]{Lu-Darvas-DiNezza-mono} used the assumption  that $\xi$ has a small unbounded locus, and explain how to remove this restriction from them. 

$1.$ Theorem 4.10 in  \cite[Section 4B]{Lu-Darvas-DiNezza-mono}: this is due to the fact that they used the integration by parts formula \cite[Theorem 1.14]{BEGZ} which requires the small unbounded locus hypothesis. Theorem \ref{th-integrabypart} frees us from this assumption and the rest of the proof of \cite[Theorem 4.10]{Lu-Darvas-DiNezza-mono} works verbatim in our setting.

$2.$ Lemma 4.17 in \cite[Section 4C]{Lu-Darvas-DiNezza-mono} and Proposition 4.30  in \cite[Section 4D]{Lu-Darvas-DiNezza-mono}: in order to  prove these results,  the authors used Lemmas 4.6 and 4.9 from Section 4A there respectively.  But these lemmas were extended to the current setting in \cite[Corollaries 3.8 and 3.9]{Lu-Darvas-DiNezza-logconcave}. Hence, Sections 4C and 4D in  \cite{Lu-Darvas-DiNezza-mono} still work in our setting.

This ends the proof.
\endproof


\bibliography{biblio_family_MA,biblio_Viet_papers}
\bibliographystyle{siam}

\bigskip

\noindent
\Addresses
\end{document}